\documentclass[12pt,a4paper]{amsart}
\usepackage[latin1]{inputenc}
\usepackage{latexsym}
\usepackage{color,graphicx,shortvrb}
\usepackage{amsmath, amssymb}
\usepackage{amsfonts}
\usepackage[colorlinks, bookmarks=true]{hyperref}

\newtheorem{theorem}{Theorem}[section]
\newtheorem{lemma}[theorem]{Lemma}
\newtheorem{proposition}[theorem]{Proposition}
\newtheorem{corollary}[theorem]{Corollary}

\theoremstyle{definition}
\newtheorem{definition}[theorem]{Definition}
\newtheorem{remark}[theorem]{Remark}

\author{J. M. Almira$^*$, Kh. F. Abu-Helaiel}
\title{A note on monomials}
\thanks{$^*$ Corresponding author}

\begin{document}


\subjclass[2010]{}


\begin{abstract} We study discontinuous solutions of the monomial equation $\frac{1}{n!}\Delta_h^nf(x)=f(h)$. In particular, we characterize the closure of their graph, $\overline{G(f)}^{\mathbb{R}^2}$, and we use the properties of these functions to present a new proof of the Darboux type theorem for polynomials and of Hamel's theorem for additive functions.
\end{abstract}

\maketitle

\markboth{J. M. Almira, Kh. F. Abu-Helaiel}{A note on monomials}

\section{Motivation}

One of the best known functional equations that exists in the literature is Fr\'{e}chet's functional equation
\begin{equation}\label{frechet}
\Delta^{n+1}_hf(x)=0 \ \ (x,h\in\mathbb{R}),                                                                                      \end{equation}
where $\Delta_h^1f(x)=f(x+h)-f(x)$ and $\Delta_{h}^{k+1}f(x)=\Delta_{h}^1\left(\Delta_{h}^kf\right)(x)$, $k=1,2,\cdots$. The solutions of this equation are named polynomials. The equation (\ref{frechet}) has been studied by many authors and all its basic properties are already known for many years. In particular, its regularity properties are well known. These can be summarized, for the case of real functions of a single real variable,  with the single statement that if $f$ is a solution of $(\ref{frechet})$, then $f$ is an ordinary polynomial of degree $\leq n$, $f(x)=a_0+a_1x+\cdots+a_nx^n$, if and only if $f$ is bounded on some set $A\subset\mathbb{R}$  with positive Lebesgue measure $|A|>0$. In particular, all measurable polynomials are ordinary polynomials. This result was firstly proved for the Cauchy functional equation by Kormes in 1926 \cite{kormes}. Later on, in 1959, the result was proved for polynomials by Ciesielski \cite{ciesielski} (see also \cite{kurepa}, \cite{haruki}). A weaker result is the so called Darboux type theorem, which claims that the polynomial $f$ is an ordinary polynomial if and only if $f_{|(a,b)}$ is  bounded for some nonempty open interval $(a,b)$ (see \cite{darboux} for the original result, which was stated for solutions of the Cauchy functional equation and \cite{almira_antonio} for a direct proof of this result with polynomials).

In \cite{almira_antonio} Fr\'{e}chet's equation was studied from a new fresh perspective. The main idea was to use the basic properties of Lagrange interpolation polynomials in one real variable. This allowed the authors to give a description of the closure of the graph $G(f)=\{(x,f(x)):x\in\mathbb{R}\}$ of any discontinuous polynomial $f$. Concretely, they proved that
\begin{equation}\label{grafo}
\overline{G(f)}^{\mathbb{R}^2}=C(l,u)=\{(x,y)\in\mathbb{R}^2: l(x)\leq y\leq u(x)\}
\end{equation}
for a certain  pair of functions  $l,u:\mathbb{R}\to\mathbb{R}\cup\{+\infty,-\infty\}$ such that
\begin{itemize}
\item[(i)] $u$ is lower semicontinuous and $l$ is upper semicontinuous.
\item[(ii)] For all $x\in \mathbb{R}$ we have that $u(x)-l(x)=+\infty$.
\item[(iii)] There exists two non-zero ordinary polynomials $p,q\in\Pi_n$ such that, $p\neq q$ and for all $x\in\mathbb{R}$, $\{x\}\times [p(x),q(x)]\subseteq C(l,u)$.                                           \end{itemize}
Clearly, this result implies the Darboux type theorem for the Fr\'{e}chet functional equation. Furthermore, it states that, for every discontinuous polynomial $f$, the set $\overline{G(f)}^{\mathbb{R}^2}$ contains an unbounded open set. This is a nice property which stands up, in a very visual form, the fact that discontinuous polynomials have wild oscillations.

In this paper we are interested on the study of  the sets $\overline{G(f)}^{\mathbb{R}^2}$  whenever $f:\mathbb{R}\to\mathbb{R}$  is a monomial. Recall that $f:\mathbb{R}\to\mathbb{R}$ is an $n$-monomial if it is a solution of the so called monomial functional equation
\begin{equation}
\label{monomials}
\frac{1}{n!}\Delta^{n}_hf(x)=f(h) \ \ (x,h\in\mathbb{R}).
\end{equation}
It is known that $f$ satisfies $(\ref{monomials})$ if and only if $f(x)=F(x,\cdots,x)$ for a certain multi-additive and symmetric function $F:\mathbb{R}^n\to\mathbb{R}$, and that $f$ is a polynomial if and only if $f(x)=\sum_{k=0}^nf_k(x)$, where $f_k(x)$ is a $k$-monomial for $k=0,1,\cdots, n$. Finally, if $f$ is an $n$-monomial, and we denote by $\Delta_{x_1x_2\cdots x_s}$ the operator given by $\Delta_{x_1x_2\cdots x_s}f(x)=\Delta_{x_1}\left(\Delta_{x_2\cdots x_s}f\right)(x)$, $s=2,3,\cdots$, then $\Delta_{x_1x_2\cdots x_n}f(x)$ does not depend on the variable $x$, and the function  $F(x_1,x_2,\cdots,x_n)= \frac{1}{n!}\Delta_{x_1x_2\cdots x_n}f(x)$  is multi-additive and symmetric (See, for example, \cite{czerwik}, \cite{kuczma}, for the proofs of these claims).

The main goal of this note is to prove that, for the case of $n$-monomials $f$, the functions $l,u$ appearing in $(\ref{grafo})$ are of the form $Ax^{n}\tau_{I}^{\varepsilon}(x)$, where $A\in\mathbb{R}$, $I\in\{ ]-\infty,0[, ]0,+\infty[, \mathbb{R},\emptyset\}$, $\epsilon\in\{+,-\}$, and
\begin{equation*}
\tau_I^\epsilon(x)=\left\{
\begin{array}{ccc}
1 &  & x\in I \\
\epsilon\cdot\infty &  & x\not\in I
\end{array}
\right.     .
\end{equation*}
Thus, in this special case, if $\overline{G(f)}^{\mathbb{R}^2}\neq \mathbb{R}^2$, then at least one of the functions $l,u$ is finite over an infinite interval (that we call the domain of finiteness of the function), and both are of the form $Ax^n$ with $A\in \mathbb{R}$ over their domains of finiteness. We also give, as a consequence, a new proof of the Darboux-type theorem for the Fr\'{e}chet functional equation. Furthermore, we also present a new proof of Hamel's theorem \cite{hamel}, \cite{sanjuan}, which claims that the graph of any discontinuous  additive function $f:\mathbb{R}\to\mathbb{R}$ is a dense subset of $\mathbb{R}^2$. We end the paper with a conjecture about the regularity properties of the functions $l,u$ for the case of general polynomials.

It is important to note that, although the Darboux type theorem is not the finest expression of the regularity property for polynomials, it is still a non-trivial result, and our proofs avoid the use of some technical results from measure theory. Furthermore, our approach has the advantage that it gives a very precise description of the closure of the graph of $f$ whenever $f$ is a discontinuous monomial.


\section{Main results}

\begin{lemma} \label{uno} If $f:\mathbb{R}\to\mathbb{R}$ is a solution of Fr\'{e}chet functional equation and $x_0,h_0\in\mathbb{R}$ then there exists a unique polynomial $p_{x_0,h_0}\in\Pi_n$ such that $f_{|x_0+h_0\mathbb{Q}}=p_{|x_0+h_0\mathbb{Q}}$.
\end{lemma}

\noindent \textbf{Proof. }  Let $f:\mathbb{R}\to\mathbb{R}$ be such that $\Delta_h^{n+1}f(x)=0$ for all $x,h\in\mathbb{R}$. Let $x_0,h_0\in\mathbb{R}$ and let $p_0(t)\in\mathbb{R}[t]$ be the polynomial of degree $\leq n$ such that $f(x_0+kh_0)=p_0(x_0+kh_0)$ for all $k\in\{0,1,\cdots,n\}$ (this polynomial exists and it is unique, thanks to Lagrange's interpolation formula). Then
\begin{eqnarray*}
0 &=&
\Delta_{h_0}^{n+1}f(x_0)=\sum_{k=0}^n\binom{n+1}{k}(-1)^{n+1-k}f(x_0+kh_0)+f(x_0+(n+1)h_0)\\
&=& \sum_{k=0}^n\binom{n+1}{k}(-1)^{n+1-k}p_0(x_0+kh_0)+f(x_0+(n+1)h_0)\\
&=& -p_0(x_0+(n+1)h_0)+f(x_0+(n+1)h_0),
\end{eqnarray*}
since  $0=\Delta_{h_0}^{n+1}p(x_0)=\sum_{k=0}^{n+1}\binom{n+1}{k}(-1)^{n+1-k}p_0(x_0+kh_0)$. This means that $f(x_0+(n+1)h_0)=p_0(x_0+(n+1)h_0)$. In particular, $p_0=q$, where $q$ denotes the polynomial of degree $\leq n$ which interpolates $f$ at the nodes $\{x_0+kh_0\}_{k=1}^{n+1}$. This argument can be repeated (forward and backward) to prove that $p_0$ interpolates $f$ at all the nodes $x_0+h_0\mathbb{Z}$.

Let $m\in\mathbb{Z}$, $m\neq 0$, and let us use the same kind of argument, taking $h_0^*=h_0/m$ instead of $h_0$. Then we get a polynomial $p_0^*$ of degree $\leq n$ such that $p_0^*$  interpolates $f$ at the nodes $x_0+\frac{h_0}{m}\mathbb{Z}$. Now, $p_0=p_0^*$ since the set
\[
\frac{h_0}{m}\mathbb{Z}\cap h_0\mathbb{Z} =  h_0\mathbb{Z}
\]
is infinite. Thus, we have proved that $p_0$ interpolates $f$ at all the points of  $$\Gamma_{x_0,h_0}:=x_0+\bigcup_{m\in\mathbb{Z}\setminus\{0\}}\frac{h_0}{m}\mathbb{Z}=x_0+h_0\mathbb{Q}.$$
{\hfill $\Box$}

\begin{lemma} \label{dos} Assume that $f:\mathbb{R}\to\mathbb{R}$ is a monomial and let $p_{x_0,h_0}$ be the polynomial described in Lemma \ref{uno}. Then
\begin{itemize}
\item[$(a)$] $p_{x_0,h_0}(t)=a_0(x_0,h_0)+a_1(x_0,h_0)t+\cdots+a_{n-1}(x_0,h_0)t^{n-1}+A_n(h_0)t^n$.
\item[$(b)$] For all $t\in h_0\mathbb{Q}$ we have that $f(t)=A_n(h_0)t^n$.
\end{itemize}
\end{lemma}

\noindent \textbf{Proof. } In principle, $p_{x_0,h_0}(t)$ is a polynomial of the form
\[
p_{x_0,h_0}(t)=a_0(x_0,h_0)+a_1(x_0,h_0)t+\cdots+a_{n-1}(x_0,h_0)t^{n-1}+A_n(x_0,h_0)t^n,
\]
so that we must prove that $A_n(x_0,h_0)=A_n(h_0)$. In other words, we must prove that $A_n(x_0,h_0)$ does not depend on the point $x_0$.

For all $m\geq 0$, and all rational number $a/b\in\mathbb{Q}$ we have that  $\Delta_{\frac{a}{b}h_0}^{m}p_{x_0,h_0}(x_0)=\Delta_{\frac{a}{b}h_0}^{m}f(x_0)$, since $f_{|x_0+h_0\mathbb{Q}}=(p_{x_0,h_0})_{|x_0+h_0\mathbb{Q}}$. Thus, a direct computation shows that
\begin{equation}\label{Ah}
A_n(x_0,h_0)\left(\frac{a}{b}h_0\right)^n = \frac{1}{n!}\Delta_{\frac{a}{b}h_0}^{n}p_{x_0,h_0}(x_0)=\frac{1}{n!}\Delta_{\frac{a}{b}h_0}^{n}f(x_0)=f(\frac{a}{b}h_0),
\end{equation}
since $f$ is an $n$-monomial. It follows that $A(x_0,h_0)$ does not depend on $x_0$ since the right hand side of $(\ref{Ah})$ does not depend on $x_0$. This proves $(a)$ and $(b)$ simultaneously. {\hfill $\Box$}

We are now in conditions to prove an ``identity principle'' for monomials:

\begin{proposition}[Identity principle] Assume that $f,g:\mathbb{R}\to\mathbb{R}$ are monomials and there exists an open interval  $I\subseteq\mathbb{R}$ such that $f(t)=g(t)$ for all $t\in I$. Then $f(t)=g(t)$ for all $t\in\mathbb{R}$.  In particular, if $f$ satisfies $\frac{1}{n!}\Delta_h^nf(x)=f(h)$ whenever $\{x,x+h,\cdots,x+nh\}\subset I$, then there exists a unique $n$-monomial $q:\mathbb{R}\to\mathbb{R}$ such that $q_{|I}=f$.
\end{proposition}

\noindent \textbf{Proof. } Let us assume, with no loss of generality, that $0\not\in I$. Let us now concentrate our attention on $f$. If we use the notation of  Lemma \ref{dos}, we know that, for all $h\in I$ and $r\in\mathbb{Q}$, $f(hr)=A_n(h)h^nr^n$. In particular, $A_n(h)=f(h)/h^n$ and
\[
f(t)=\frac{f(h)}{h^n}(t)^n \text{ for all } t\in h\mathbb{Q}.
\]
This means that the value of $f$ at $h\in I$ uniquely determines the value of $f$ at $t$ for all $t\in h\mathbb{Q}$. Now, $\bigcup_{h\in I}h\mathbb{Q}=\mathbb{R}$, since $I$ is open. This implies that $f$ is uniquely defined over $\mathbb{R}$. {\hfill $\Box$}


\begin{lemma} \label{tres} Assume that $f:\mathbb{R}\to\mathbb{R}$ is a discontinuous $n$-monomial and let $A_n(h)$ denote the constant appearing in part $(b)$ of Lemma \ref{dos}. Then
\begin{itemize}
\item[$(a)$] If $h_0,h_1\in\mathbb{R}$ satisfy  $A_n(h_0)\neq A_n(h_1)$, the space between the monomials $A_n(h_0)t^n$ and $A_n(h_1)t^n$ is a subset of $\overline{G(f)}^{\mathbb{R}^2}$.

\item[$(b)$] $\sup_{h\in\mathbb{R}}|A_n(h)|=+\infty$.
\end{itemize}
\end{lemma}

\noindent \textbf{Proof. }Let us prove $(a)$. Take $h_0,h_1\in\mathbb{R}$ such that $A_n(h_0)\neq A_n(h_1)$. These values obviously exists since $f$ is assumed to be discontinuous and $p_{0,h}(t)=A_n(h)t^n$. Take $x_0\in h_0\mathbb{Q}\setminus \{0\}$ and $x_1\in h_1\mathbb{Q}\setminus \{0\}$. Take $i\in\{0,1,\cdots,n\}$, $r\in\mathbb{Q}$ and set $x_{i,r}=(1-i)x_0+i r x_1$. Obviously, $x_{i,r}\in (1-i)x_0+ x_1\mathbb{Q}$ for all $r\in\mathbb{Q}$. Thus, if we define $p_0=p_{0,h_0}$, and $p_i=p_{(1-i)x_0,x_1}$, for $i=1,2,\cdots,n$, then $p_i(x_{i,r})=f(x_{i,r})$ for all $r\in\mathbb{Q}$ and all $i\in\{0,1,\cdots,n\}$. Furthermore, $p_0(t)=p_{0,h_0}(t)=A_n(h_0)t^n$ and  $p_1(t)=p_{0,x_1}(t)=p_{0,h_1}(t)=A_n(h_1)t^n$.
Let us now consider, for each $r\in\mathbb{Q}$ the unique polynomial $q_{r}(t)$ of degree $\leq n$ which interpolates to $f$ at the nodes $\{x_{i,r}\}_{i=0}^n$.

If we take into account that $x_{i,r}$ can be expressed as $x_{i,r}=x_0+i(rx_1-x_0)$,  we conclude that $q_r=p_{x_0,(rx_1-x_0)}$, which means that the graph of $q_r$ is a subset of $\overline{G(f)}^{\mathbb{R}^2}$. Let us now use that $\mathbb{Q}$ is a dense subset of $\mathbb{R}$ to force $r\in\mathbb{Q}$ to tend to $\frac{x_0}{x_1}$. This has the effect that the set $\{x_{i,r}\}$ collapses to the point $\{x_0\}$, since $\lim_{r\in\mathbb{Q},r\to\frac{x_0}{x_1}}x_{i,r}=x_0$ , $i=0,1,\cdots,n$. This should force the graph of $q_{r}$ to blow up when $r\to\frac{x_0}{x_1}$.

Concretely, we have that
\begin{eqnarray*}
\lim_{r\in\mathbb{Q},r\to\frac{x_0}{x_1}}q_r(x_{i,r})&=& \lim_{r\in\mathbb{Q},r\to\frac{x_0}{x_1}}f(x_{i,r})\\
&=& \lim_{r\in\mathbb{Q},r\to\frac{x_0}{x_1}}p_i(x_{i,r}) \\
&=& p_i(x_0), \ \ i=0,1,\cdots,n.
\end{eqnarray*}
Thus, if we set $q(x)=L[\{i\}_{i=0}^n,\{p_i(x_0)\}_{i=0}^n](x)$ and  $\widetilde{q}_r(x)=L[\{i\}_{i=0}^n,\{q_r(x_{i,r})\}_{i=0}^n](x)$, where $L[\{x_k\}_{k=0}^n,\{y_k\}_{k=0}^n](x)$ denotes the Lagrange interpolation polynomial associated to the nodes $\{x_k\}_{k=0}^n$ and the values $\{y_k\}_{k=0}^n$, then $\widetilde{q}_r$ converges (for $r\in\mathbb{Q},r\to\frac{x_0}{x_1}$)  uniformly on compact subsets of the real line, to the polynomial $q$.
Note that
\[
i=\frac{x_{i,r}-x_0}{rx_1-x_0}, \text{ for } i=0,1,\cdots,n,
\]
so that:
\begin{equation*}
q_r(x_{i,r})=\widetilde{q}_r(i)= \widetilde{q}_r\left(\frac{x_{i,r}-x_0}{rx_1-x_0}\right), \text{ for } i=0,1,\cdots, n.\end{equation*}
and  \begin{equation}  \label{qr}
q_r(x)= \widetilde{q}_r\left(\frac{x-x_0}{rx_1-x_0}\right).
\end{equation}

Obviously, $J_0=[A_n(h_0)x_0^n,A_n(h_1)x_0^n]=[p_0(x_0),p_1(x_0)]\subseteq q(\mathbb{R})$. since $q$ is continuous, $q(0)=p_0(x_0)$, and $q(1)=p_1(x_0)$. We want to show that $\{x_0\}\times J_0$ is a subset of $\overline{G(f)}^{\mathbb{R}^2}$. In fact, we prove much more than this, since we demonstrate that $\{x_0\}\times q(\mathbb{R}) \subseteq \overline{G(f)}^{\mathbb{R}^2}$.

Take $c\in q(\mathbb{R})$ and let $\varepsilon>0$ be a positive number. We can find $s\in\mathbb{Q}$ such that $|q(s)-c|<\varepsilon/2$, since $\mathbb{Q}$ is dense in $\mathbb{R}$. Let $x_{s,r}=s(rx_1-x_0)+x_0$. Then $q_r(x_{s,r})=f(x_{s,r})$. Furthermore, formula  $(\ref{qr})$ tell us that
\[
q_r(x_{s,r})= \widetilde{q}_r\left(\frac{x_{s,r}-x_0}{rx_1-x_0}\right)= \widetilde{q}_r(s).
\]
Now, $\lim_{r\in\mathbb{Q},r\to\frac{x_0}{x_1}}\|\widetilde{q}_r-q\|_{[0,s+1]}=0$, so that, for $|r- \frac{x_0}{x_1}|$ small enough,
\[
|f(x_{s,r})-c|=|\widetilde{q}_r(s)-c|\leq |\widetilde{q}_r(s)-q(s)|+ |q(s)-c|<\varepsilon.
\]
Hence $(x_0,c)\in \overline{G(f)}^{\mathbb{R}^2}$, which is what we wanted to prove. This demonstrates that the space between the monomials $A_n(h_0)t^n$ and $A_n(h_1)t^n$ is contained into $\overline{G(f)}^{\mathbb{R}^2}$, since $x_0\in h_0\mathbb{Q}\setminus\{0\}$ was arbitrary.

To prove part $(b)$ of the lemma, we use that $q(\mathbb{R})$ is unbounded, since $q$ is a non-constant polynomial. This implies that $f$ is locally unbounded, so that
\[
\sup_{h\neq 0}|A_n(h)|\geq \sup_{h\in [1,2]}|A_n(h)|= \sup_{h\in [1,2]}\frac{|f(h)|}{|h^n|}=+\infty.
\]
{\hfill $\Box$}

\begin{corollary}[Hamel's theorem] If $f:\mathbb{R}\to\mathbb{R}$ is additive and discontinuous, then $\overline{G(f)}^{\mathbb{R}^2}=\mathbb{R}^2$.
\end{corollary}

\noindent \textbf{Proof. } Obviously, $f$ is additive if and only if it is a $1$-monomial (i.e, $\Delta_hf(x)=f(h)$). Thus, given $f$ a discontinuous $1$-monomial, let $A_1(h)$ be the constant appearing in part $(b)$ of Lemma \ref{dos}. It follows, from part $(a)$ of Lemma \ref{tres}, that we only need to prove that,
\[
-\infty=\inf_{h\neq 0}A_1(h)<\sup_{h\neq 0}A_1(h)=+\infty.
 \]

Our assumption about the discontinuity of $f$ implies the existence of $h,s\in\mathbb{R}\setminus\{0\}$ such that $A_1(h)\neq A_1(s)$. What is more: we can assume $h,s>0$ with no lost of generality. Obviously,
\begin{eqnarray*}
A_1(h-s) &=&  \frac{f(h-s)}{h-s} = \frac{f(h)}{h-s}-\frac{f(s)}{h-s}\\
&=&  \frac{f(h)}{h}\frac{h}{h-s}-\frac{f(s)}{s}\frac{s}{h-s} \\
&=& A_1(h)\frac{h}{h-s}-A_1(s)\frac{s}{h-s}\\
&=& A_1(h)+ (A_1(h)- A_1(s))\frac{s}{h-s}.
\end{eqnarray*}
Let us now consider $A_1(h-rs)$ as a function of $r\in\mathbb{Q}$. Then
\[
A_1(h-rs)=A_1(h)+ (A_1(h)- A_1(s))\frac{rs}{h-rs} = A_1(h)+ (A_1(h)- A_1(s))\frac{1}{\frac{h}{rs}-1},
\]
since $A_1(rs)=A_1(s)$ for all $s\in\mathbb{R}\setminus\{0\}$ and all $r\in\mathbb{Q}\setminus\{0\}$. Hence,
\[
\lim_{r\to (h/s)^+,r\in\mathbb{Q}}A_1(h-rs)= (-\infty)\cdot \mathbf{sign}(A_1(h)- A_1(s))
\]
and
\[
\lim_{r\to (h/s)^-,r\in\mathbb{Q}}A_1(h-rs)= (+\infty)\cdot \mathbf{sign}(A_1(h)- A_1(s)).
\]
This ends the proof. {\hfill $\Box$}






The following result should be known to experts, but we include it here for the sake of completeness:

\begin{lemma}\label{cuatro} Assume that $f_i:\mathbb{R}\to\mathbb{R}$ is an $n_i$-monomial $(i=1,2)$. Then $f(x)=f_1(x)f_2(x)$ is an $(n_1+n_2)$-monomial. Furthermore, if $f_1,f_2:\mathbb{R}\to\mathbb{R}$ are $n$-monomials, then $f_1+f_2$ is also an $n$-monomial.
\end{lemma}
\noindent \textbf{Proof. } We prove the first claim, since the second one is obvious. By hypothesis, there are two symmetric multi-additive functions $F_i:\mathbb{R}^{n_i}\to \mathbb{R}$ such that $f_i(x)=F_i(x,x,\cdots,x)$ $(i=1,2)$. Set
\[
F(x_1,\cdots,x_{n_1+n_2})=\frac{1}{(n_1+n_2)!}\sum_{\sigma\in \mathcal{S}_{n_1+n_2}}F_1(x_{\sigma(1)},\cdots,x_{\sigma(n_1)})F_2(x_{\sigma(n_1+1)}, \cdots,x_{\sigma(n_1+n_2)}),
\]
where $\mathcal{S}_N$ denotes the group of permutations of the set $\{1,2,\cdots,N\}$. Then $F$ is symmetric and  $(n_1+n_2)$-additive, and
\[
f_1(x)f_2(x)=F(x,\cdots,x).
\]
{\hfill $\Box$}

\begin{theorem} \label{teoppal} Assume that $f$ is a discontinuous $n$-monomial and let $\Gamma=\overline{G(f)}^{\mathbb{R}^2}$. Let $\alpha=\sup_{h\in\mathbb{R}^*}A_n(h)$ and  $\beta=\inf_{h\in\mathbb{R}^*}A_n(h)$.  Then:
\begin{itemize}
\item[$(a)$] If $\alpha=+\infty$ and $\beta=-\infty$, then $\Gamma =\mathbb{R}^2$.
\item[$(b)$] If $\alpha=+\infty$ and $\beta\in\mathbb{R}$, then $\Gamma =\{(x,y):y\geq \beta x^{n}\}$ if $n=2k$ is an even number, and $\Gamma =\{(x,y):x\leq 0 \text{ and } y\leq \beta x^{n}\}\cup\{(x,y):x\geq 0\text{ and } y\geq \beta x^{n}\}$ if $n=2k+1$ is an odd number.
    In particular, if $\beta=0$, we get the half space $\Gamma =\{(x,y):y\geq 0\}$ for $n=2k$ and the union of the first and third quadrants $\Gamma =\{(x,y):xy\geq 0\}$, for $n=2k+1$.
\item[$(c)$] If $\alpha\in\mathbb{R}$ and $\beta=-\infty$, then $\Gamma =\{(x,y):y\leq \beta x^{n}\}$ if $n=2k$ is an even number, and $\Gamma =\{(x,y):x\leq 0 \text{ and } y\geq \beta x^{n}\}\cup\{(x,y):x\geq 0\text{ and } y\leq \beta x^{n}\}$ if $n=2k+1$ is an odd number.  In particular, if $\alpha=0$, we get the half space $\Gamma =\{(x,y):y\leq 0\}$ for $n=2k$ and the union of the second and fourth quadrants $\Gamma =\{(x,y):xy\leq 0\}$, for $n=2k+1$.
 \end{itemize}
Furthermore, for all $n\geq 2$ there are examples of discontinuous $n$-monomials $f$ verifying each one of the claims $(a),(b),(c)$ above.
\end{theorem}
\noindent \textbf{Proof.} Part $(b)$ of Lemma \ref{tres} implies that $\alpha=+\infty$ or $\beta=-\infty$.  Hence the three cases $(a)$, $(b)$ and $(c)$  considered in the statement of the theorem exhaust all possibilities. Let us consider each case separately:

\medskip

\noindent \textbf{Case (a): $\alpha=+\infty$ and $\beta=-\infty$. } By hypothesis, there exists two sequences of real numbers $\{\tau_k\}$, $\{\eta_k\}$ such that $A_n(\tau_k)>A_n(\eta_k)$ for all $k\in\mathbb{N}$, $\lim_{k\to\infty}A_n(\tau_k)=+\infty$, and $\lim_{k\to\infty}A_n(\eta_k)=-\infty$. Now, part $(a)$ of Lemma \ref{tres} implies that, for each $k\in\mathbb{N}$,  $\overline{G(f)}^{\mathbb{R}^2}$ contains the space between the monomials $A_n(\tau_k)t^n$ and  $A_n(\eta_k)t^n$, which obviously implies that  $\overline{G(f)}^{\mathbb{R}^2}= \mathbb{R}^2$.

\medskip

\noindent \textbf{Case (b): $\alpha=+\infty$ and $\beta\in\mathbb{R}$. } In this case, there exists  two sequences of real numbers $\{\tau_k\}$, $\{\eta_k\}$ such that $A_n(\tau_k)>A_n(\eta_k)$ for all $k\in\mathbb{N}$, $\lim_{k\to\infty}A_n(\tau_k)=+\infty$, and $\lim_{k\to\infty}A_n(\eta_k)=\beta$. Hence, we can use again part $(a)$ of Lemma \ref{tres} to prove that $\Gamma=\overline{G(f)}^{\mathbb{R}^2}$ contains the set $\Sigma_1 =\{(x,y):y\geq \beta x^{n}\}$ if $n=2k$ is an even number, or it contains the set $\Sigma_2 =\{(x,y):x\leq 0 \text{ and } y\leq \beta x^{n}\}\cup\{(x,y):x\geq 0\text{ and } y\geq \beta x^{n}\}$ if $n=2k+1$ is an odd number. We must prove, in both cases, that $\Gamma$ does not contain any other point. Thus, let us assume that $n\in 2\mathbb{N}$ and $(x_0,y_0)\in\Gamma\setminus \Sigma_1$. Then $y_0<\beta x_0^n$. We can assume, without loss of generality, that $(x_0,y_0)\in G(f)$ since $\mathbb{R}^2\setminus \Sigma_1$ is an open set. Then $A_n(x_0)=\frac{f(x_0)}{x_0^n}<\beta$, which is impossible. This proves that $\Gamma=\Sigma_1$. If $n\in 2\mathbb{N}+1$, the arguments are similar.

\noindent \textbf{Case (c): $\alpha\in\mathbb{R}$ and $\beta=-\infty$. }  This case admits a proof completely analogous to the proof given for the case (b) above.
\medskip

Let us now show that all the situations considered in cases $(a)$, $(b)$ and $(c)$ above, are supplied by concrete simple examples. For $n=1$ we always have that $\overline{G(f)}^{\mathbb{R}^2} =\mathbb{R}^2$ by Hamel's theorem. Thus, we will assume in all what follows that $n\geq 2$.

To proceed with this part of the proof, we need first to introduce a concrete example of discontinuous additive function. Let $\Upsilon=\{s_i\}_{i\in I}$ be an algebraic basis of $\mathbb{R}$ as a $\mathbb{Q}$-vector space. Then any map $\phi:\Upsilon\to\mathbb{R}$ can be uniquely extended as a $\mathbb{Q}$-linear map to a real function of one real variable $\Phi:\mathbb{R}\to\mathbb{R}$. We denote by $\mathcal{L}:\mathbb{R}\to\mathbb{R}$ the unique $\mathbb{Q}$-linear map satisfying $\mathcal{L}(s_i)=1$ for all $i\in I$. Obviously, $\mathcal{L}(x)$ is a  discontinuous $1$-monomial, so that Lemma \ref{cuatro} implies that  $ax^n+bx^k\mathcal{L}(x)^{n-k}$ is an $n$-monomial for all $a,b\in\mathbb{R}$ and $k=0,1,\cdots, n-1$.

Let $n\in\mathbb{N}$ be fixed. Then $f(x)=x^{n-1}\mathcal{L}(x)$ is an $n$-monomial such that $\overline{G(f)}^{\mathbb{R}^2}= \mathbb{R}^2$. To prove this, we only need to check that $\alpha=\sup_{h\in\mathbb{R}^*}A_n(h)=+\infty$ and  $\beta=\inf_{h\in\mathbb{R}^*}A_n(h)=-\infty$. Now, $\mathcal{L}(x)$  has a dense graph in the plane, so that $\sup_{x\in [1,2]}\mathcal{L}(x)=+\infty$ and $\inf_{x\in [1,2]}\mathcal{L}(x)=-\infty$. Hence $\alpha\geq \sup_{h\in [1,2]}A_n(h)= \sup_{h\in [1,2]}\frac{\mathcal{L}(h)h^{n-1}}{h^n}=\sup_{h\in [1,2]}\frac{\mathcal{L}(h)}{h}=+\infty$ and, analogously,  $\beta \leq \inf_{h\in [1,2]}A_n(h)= \inf_{h\in [1,2]}\frac{\mathcal{L}(h)}{h}=-\infty$. This proves the existence of $n$-monomials satisfying case (a). For case (b), we may choose $f(x)= \beta x^n+\mathcal{L}(x)^n$ whenever $n=2k$ is an even number, and $f(x)= \beta x^n+x\mathcal{L}(x)^{n-1}$ whenever $n=2k+1$ is an odd number. Case (c) is covered by the examples:  $f(x)= \alpha x^n-\mathcal{L}(x)^n$ whenever $n=2k$ is an even number, and $f(x)= \alpha x^n-x\mathcal{L}(x)^{n-1}$ whenever $n=2k+1$ is an odd number.

{\hfill $\Box$}

\begin{remark} Theorem \ref{teoppal} may be reformulated, in terms of the functions $l,u$ appearing in $(\ref{grafo})$, as follows: Let $f$ be a discontinuous $n$-monomial  and let $\Gamma=\overline{G(f)}^{\mathbb{R}^2}=C(l,u)$ be the representation of the closure of the graph of $f$ given by  $(\ref{grafo})$.   Let $\alpha=\sup_{h\in\mathbb{R}^*}A_n(h)$ and  $\beta=\inf_{h\in\mathbb{R}^*}A_n(h)$.  Then:
\begin{itemize}
\item[$(a)$] If $\alpha=+\infty$ and $\beta=-\infty$, then $l(x)=-\infty$ and $u(x)=+\infty$ for all $x\in\mathbb{R}$.
\item[$(b)$] If $\alpha=+\infty$, $\beta\in\mathbb{R}$, and $n=2k$ is an even number, then $l(x)=\beta x^n$ and $u(x)=+\infty$ for all $x\in\mathbb{R}$.
\item[$(c)$]If $\alpha=+\infty$, $\beta\in\mathbb{R}$, and $n=2k+1$ is an odd number, then
\begin{equation*}
l(x)=\left\{
\begin{array}{ccc}
-\infty &  & x\leq 0 \\
\beta x^n &  & x>0
\end{array}
\right.     \text{ and }
u(x)=\left\{
\begin{array}{ccc}
\beta x^n  &  & x< 0 \\
+\infty &  & x\geq 0
\end{array}
\right.  .
\end{equation*}
\item[$(d)$] If $\alpha\in\mathbb{R}$, $\beta=-\infty$, and $n=2k$ is an even number, then $l(x)=-\infty$ and $u(x)=\alpha x^n$ for all $x\in\mathbb{R}$.
\item[$(e)$]If $\alpha\in\mathbb{R}$, $\beta=-\infty$, and $n=2k+1$ is an odd number, then
\begin{equation*}
l(x)=\left\{
\begin{array}{ccc}
\alpha x^n &  & x< 0 \\
-\infty &  & x\geq 0
\end{array}
\right.     \text{ and }
u(x)=\left\{
\begin{array}{ccc}
+\infty   &  & x\leq  0 \\
\alpha x^n &  & x> 0
\end{array}
\right.  .
\end{equation*}
   \end{itemize}
Furthermore, for all $n\geq 2$ there are examples of discontinuous $n$-monomials $f$ verifying each one of the claims $(a),(b),(c), (d), (e)$ above.
\end{remark}

\begin{corollary} If $f:\mathbb{R}\to\mathbb{R}$ is a discontinuous monomial of even degree, then $F(x,h)=\Delta_hf(x)$ satisfies $\overline{G(F)}^{\mathbb{R}^3}=\mathbb{R}^3$.
\end{corollary}

\noindent \textbf{Proof. } It follows from Theorem \ref{teoppal} that $\Gamma= \overline{G(f)}^{\mathbb{R}^2}$ satisfies $\Gamma=\mathbb{R}^2$, or $\Gamma=\{(x,y)\in\mathbb{R}^2: y\geq \beta x^n\}$, or  $\Gamma =\{(x,y)\in\mathbb{R}^2: y\leq \alpha x^n\}$.  Let us assume that $\Gamma=\{(x,y)\in\mathbb{R}^2: y\geq \beta x^n\}$ (the other cases admit similar proofs).

Let $(x,h,\lambda)\in\mathbb{R}^3$ and take $\varepsilon >0$. Obviously, there exists a constant $M\in\mathbb{R}$ such that $[x,x+h]\times [M, +\infty) \subset \Gamma$. Take $a,b\in [M,+\infty)$ such that  $\lambda=b-a$. By construction, there exists two real numbers $x^*,h^*$ such that $\max\{|x-x^*|,|h-h^*|,|f(x^*)-a|,|f(x^*+h^*)-b|\}<\frac{\varepsilon}{2}$. Hence
\begin{eqnarray*}
\|(x^*,h^*,F(x^*,h^*))-(x,h,\lambda)\|_{\infty} &=& \|(x^*-x,h^*-h,f(x^*+h^*)-f(x^*)-\lambda)\|_{\infty}\\
&=& \max\{|x-x^*|,|h-h^*|,|f(x^*+h^*)-b-(f(x^*)-a)|\}\\
&\leq& \max\{|x-x^*|,|h-h^*|,|f(x^*+h^*)-b|+|f(x^*)-a|\} \leq \varepsilon.
\end{eqnarray*}
This proves that $(x,h,\lambda)\in\overline{G(F)}^{\mathbb{R}^3}$. Hence $\overline{G(F)}^{\mathbb{R}^3}$, since  $(x,h,\lambda)$ was arbitrary. {\hfill $\Box$}

\begin{corollary}[Darboux type Theorem for Fr\'{e}chet's functional equation]  Let $f$ be a polynomial, and let $f(x)=f_0+f_1(x)+\cdots+f_N(x)$ be the decomposition of $f$ as sum of monomials. If $f_N(x)$ is a discontinuous monomial then $f$ is locally unbounded (i.e., for all $a<b$ we have that $f([a,b])$ is an unbounded subset of $\mathbb{R}$).  Consequently, the polynomial $f$ is bounded over some open interval if and only if is an ordinary polynomial, $f(x)=a_0+a_1x+\cdots+a_Nx^N$.
\end{corollary}

\noindent \textbf{Proof. }Let us prove the first claim. Assume, on the contrary, that $f_{|[a,b]}$ is bounded for a certain interval $[a,b]$ with $a<b$. Taking differences, we have that
\[
\frac{1}{N!}\Delta_h^{N}f(a)=f_N(h),
\]
so that $f_{N}$ should be bounded over the set $[0,\frac{b-a}{N}]$, which, thanks to Theorem \ref{teoppal}, would imply that $f_N$ is a continuous monomial, $f_N(x)=A_nx^N$.

Let us now assume that that $f_{|[a,b]}$ is bounded. Let $j_0=\max\{j\in\{0,1,\cdots,N\}:f_j \text{ is discontinuous}\}$. Then $g(x)=f(x)-(f_{j_0+1}(x)+\cdots+f_N(x))$ is bounded on $[a,b]$ and $g(x)=f_0+f_1(x)+\cdots+f_{j_0}(x)$ is the decomposition of $g$ as sum of monomials. The result follows just applying the first part of this corollary.
{\hfill $\Box$}


\begin{remark} It is interesting to observe that, in general, the closure of the graph of the sum of two discontinuous $n_k$-monomials $f_k$,  $(k=1,2)$, with $n_1\neq n_2$, is not uniquely determined by the closures of the graphs of these monomials. For example, we can set $f_1(x)=x^2+\mathcal{L}(x)^2$, $f_2(x)=x^2+2\mathcal{L}(x)^2$, $f_3(x)=x\mathcal{L}(x)^2$, $f_4(x)=f_1(x)+f_3(x)=x^2+(1+x)\mathcal{L}(x)^2$, and $f_5(x)=f_2(x)+f_3(x)=x^2+(2+x)\mathcal{L}(x)^2$. Then $f_1$ and $f_2$  are $2$-monomials, $f_3$ is a $3$-monomial, and the sets $\Gamma_k=\overline{G(f_k)}^{\mathbb{R}^2}$, $k=1,2,3,4,5$ are given by
\begin{eqnarray*}
\Gamma_1 = \Gamma_2 &=&  \{(x,y)\in\mathbb{R}^2: y\geq x^2\}\\
\Gamma_3 &=&  \{(x,y)\in\mathbb{R}^2: xy\geq 0\} \\
\Gamma_4 &=&  \{(x,y)\in\mathbb{R}^2: x\leq -1 \text{ and } y\leq x^2\} \cup \{(x,y)\in\mathbb{R}^2: x\geq -1 \text{ and } y\geq x^2\}\\
\Gamma_5 &=&  \{(x,y)\in\mathbb{R}^2: x\leq -2 \text{ and } y\leq x^2\} \cup \{(x,y)\in\mathbb{R}^2: x\geq -2 \text{ and } y\geq x^2\}.
\end{eqnarray*}
This makes very difficult to say something about the closure of the graph of a general polynomial just taking into account what happens for monomials. In any case, we conjecture that the functions $l,u$ appearing in $(\ref{grafo})$ for $f$ any general polynomial are, in their domains of finiteness, continuous splines.
\end{remark}

 \bibliographystyle{amsplain}

\begin{thebibliography}{99}
\bibitem{almira_antonio} \textbf{J. M. Almira, A. J.  L\'{o}pez-Moreno, } On solutions of the Fr\'{e}chet functional equation, J. Math. Anal. Appl.  332 (2007), 1119--1133.

\bibitem{almira_frechet_padico} \textbf{J. M. Almira, } A note on classical and $p$-adic Fr\'{e}chet functional equation with restrictions, Results in Math., 2011, in press, DOI 10.1007/s00025-011-0223-9.

\bibitem{cauchy}  \textbf{A. L. Cauchy, } \textit{Cours d'analyse de l'Ecole
Polytechnique, 1.} Analyse Alg\'{e}brique, V. Par\'{\i}s, 1821. [Oeuvres (2) \textbf{3%
}, Par\'{\i}s, 1897].

\bibitem{ciesielski} \textbf{Z. Ciesielski, } Some properties of convex functions of higher order, Ann. Pol. Math. \textbf{7} (1959) 1-7.
\bibitem{czerwik} \textbf{S. Czerwik, } \textit{Functional equations and inequalities in several variables,} World Scientific, 2002.
\bibitem{darboux} \textbf{G. Darboux, } Memoire sur les fonctions
discontinues, Ann. Sci. Scuola. Sup. \textbf{4} (1875) 57-112.

\bibitem{frechet} \textbf{M. Fr\'{e}chet, }   Une definition fonctionelle des polynomes,
Nouv. Ann. \textbf{9} (1909), 145-162.

\bibitem{ger1}  \textbf{R. Ger}, On some properties of polynomial functions,
Ann. Pol. Math. \textbf{25 }(1971) 195-203.
\bibitem{ger} \textbf{R. Ger, } On extensions of polynomial functions, Results in
Mathematics \textbf{26} (1994), 281-289.

\bibitem{hamel}  \textbf{G. Hamel}, Einer basis aller Zahlen und die
unstetigen L-sungen der Funktionalgleichung $f(x+y)=f(x)+f(y)$,
Math. Ann. \textbf{60} (1905) 459-472.

\bibitem{haruki} \textbf{S. Haruki, } On the theorem of S. Kakutani-M. Nagumo and J.L. Walsh for the mean value property of harmonic and complex polynomials, Pacific J. Math. \textbf{94} (1) (1981) 113-123

\bibitem{kormes} \textbf{M. Kormes, } On the functional equation $f(x+y)=f(x)+f(y)$, Bulletin of the Amer. Math. Soc. \textbf{32} (1926) 689-693.
\bibitem{kuczma} \textbf{M. Kuczma}, \textit{An introduction to the theory of functional equations and inequalities, } Polish Scientific Publ. and Silesian University Press, 1985.

\bibitem{kuczma1} \textbf{M. Kuczma}, On measurable functions with
vanishing differences, Ann. Math. Sil. \textbf{6} (1992) 42-60.
\bibitem{kurepa} \textbf{S. Kurepa, } A property of a set  of positive measure and its application, J. Math. Soc. Japan \textbf{13} (1) (1961) 13-19.

\bibitem{mckiernan}  \textbf{M. A. Mckiernan, }On vanishing n-th ordered
differences and Hamel bases, Ann. Pol. Math. \textbf{19} (1967) 331-336.

\bibitem{popoviciu} \textbf{T. Popoviciu, } Sur quelques propriétés des fonctions d'une ou deux variables reélles, Mathematica (Cluj) \textbf{8} (1934) 1-85.

\bibitem{sanjuan} \textbf{R. San Juan, } Una aplicación de las aproximaciones diofánticas a la ecuación funcional $f(x_1+x_2)=f(x_1)+f(x_2)$, Publicaciones del Inst. Matemático de la Universidad Nacional del Litoral \textbf{6} (1946) 221-224.

\end{thebibliography}


\bigskip

\footnotesize{J. M. Almira

Departamento de Matem\'{a}ticas. Universidad de Ja\'{e}n.

E.P.S. Linares,  C/Alfonso X el Sabio, 28

23700 Linares (Ja\'{e}n) Spain

Email: jmalmira@ujaen.es }



\vspace{1cm}

\footnotesize{Kh. F. Abu-Helaiel

Departamento de Estadística e Investigación Operativa. Universidad de Ja\'{e}n.

Campus de Las Lagunillas

23071 Ja\'{e}n, Spain

Email: kabu@ujaen.es
}



\end{document}